\documentclass[amstex,12pt,russian,amssymb]{article}

\usepackage{mathtext}
\usepackage[cp1251]{inputenc}
\usepackage[T2A]{fontenc}
\usepackage[russian]{babel}
\usepackage[dvips]{graphicx}
\usepackage{amsmath}
\usepackage{amssymb}
\usepackage{amsxtra}
\usepackage{latexsym}
\usepackage{ifthen}

\begin{document}

\newcounter{lemma}
\newcommand{\lemma}{\par \refstepcounter{lemma}%
{\bf Лемма \arabic{lemma}.}}

\newcounter{corollary}
\newcommand{\corollary}{\par \refstepcounter{corollary}%
{\bf Следствие \arabic{corollary}.}}

\newcounter{remark}
\newcommand{\remark}{\par \refstepcounter{remark}%
{\bf Замечание \arabic{remark}.}}

\newcounter{theorem}
\newcommand{\theorem}{\par \refstepcounter{theorem}%
{\bf Теорема \arabic{theorem}.}}

\newcounter{proposition}
\newcommand{\proposition}{\par \refstepcounter{proposition}%
{\bf Предложение \arabic{proposition}.}}

\newcommand{\proof}{{\it Доказательство.\,\,}}
\renewcommand{\refname}{\centerline{\bf Список литературы}}

{\bf Р.Р.~Салимов} (Институт математики НАН Украины),

{\bf Е.А.~Севостьянов} (Житомирский государственный университет
им.~И.~Франко)

\medskip
{\bf Р.Р.~Салімов} (Інститут математики НАН України),

{\bf  Є.О.~Севостьянов} (Житомирський державний університет
ім.~І.~Франко)

\medskip
{\bf R.R.~Salimov} (Institute of Mathematics of NAS of Ukraine),

{\bf E.A.~Sevost'yanov} (Zhitomir Ivan Franko State University)

\medskip
{\bf О равностепенной непрерывности обратных отображений в замыкании
области}

\medskip
Изучается поведение обратных гомеоморфизмов для класса отображений,
в котором выполнены верхние оценки модуля семейств кривых. В
терминах простых концов пространственных областей доказано, что
семейства таких гомеоморфизмов равностепенно непрерывны (нормальны)
в замыкании заданной области.

\medskip
{\bf Про одностайну неперервність обернених відображень в замиканні
області}

\medskip
Вивчається локальна поведінка обернених гомеоморфізмів для класу
відображень, в котрому виконано верхні оцінки модуля сімей кривих. В
термінах простих кінців просторових областей доведено, що сім'ї
таких гомеоморфизмів одностайно неперервні (нормальні) в замиканні
заданої області.

\medskip
{\bf On equicontinuity of inverse mappings in a closure of a domain}

\medskip
For some class of mappings satisfying upper modular estimates with
respect to families of curves, a behavior of the corresponding
inverse mappings is investigated. In the terms of prime ends, it is
proved that, families of such homeomorphisms are equicontinuous
(normal) in the closure of a given domain.

\newpage
{\bf 1. Введение.} Основные определения и обозначения,
использующиеся ниже, могут быть найдены в монографии \cite{MRSY}
либо статье \cite{KR}.

В нашей сравнительно недавней публикации \cite{Sev$_1$} установлено
свойство равностепенной  непрерывности для отображений, обратные к
которым являются так называемыми $Q$-гомеоморфизмами -- наиболее
простейшими обобщениями квазиконформных отображений по О. Мартио
(см. главу 4 в \cite{MRSY}). Отметим, что речь идёт здесь о
равностепенной непрерывности с локально связными границами, однако,
для областей с более общими типами границ данный вопрос до сих пор
не исследован. В данной работе мы несколько усилим упомянутые
результаты, рассматривая более широкие типы областей, для которых
указанные утверждения всё ещё имеют место. Ниже пойдёт речь о
равностепенной непрерывности упомянутых отображений в терминах
простых концов, так как даже их непрерывное продолжение на границу в
поточечном смысле, вообще говоря, не гарантируется.

Напомним некоторые определения, а также приведём формулировку
основных результатов работы. Всюду далее $D$ -- область в ${\Bbb
R}^n,$ $n\ge 2,$ $m$ -- мера Лебега в ${\Bbb R}^n.$

Следующие определения могут быть найдены в работе \cite{KR}. Пусть
$\omega$ -- открытое множество в ${\Bbb R}^k$, $k=1,\ldots,n-1$.
Непрерывное отображение $\sigma:\omega\rightarrow{\Bbb R}^n$
называется {\it $k$-мерной поверхностью} в ${\Bbb R}^n$. {\it
Поверхностью} будет называться произвольная $(n-1)$-мерная
поверхность $\sigma$ в ${\Bbb R}^n.$ Поверхность
$\sigma:\omega\rightarrow D$ называется {\it жордановой
поверхностью} в $D$, если $\sigma(z_1)\ne\sigma(z_2)$ при $z_1\ne
z_2$. Далее мы иногда будем использовать $\sigma$ для обозначения
всего образа $\sigma(\omega)\subset {\Bbb R}^n$ при отображении
$\sigma$, $\overline{\sigma}$ вместо $\overline{\sigma(\omega)}$ в
${\Bbb R}^n$ и $\partial\sigma$ вместо
$\overline{\sigma(\omega)}\setminus\sigma(\omega)$. Жорданова
поверхность $\sigma$ в $D$ называется {\it разрезом} области $D$,
если $\sigma$ разделяет $D$, т.\,е. $D\setminus \sigma$ имеет ровно
две компоненты, $\partial\sigma\cap D=\varnothing$ и
$\partial\sigma\cap\partial D\ne\varnothing$.

Последовательность $\sigma_1,\sigma_2,\ldots,\sigma_m,\ldots$
разрезов области $D$ называется {\it цепью}, если:

\medskip

(i) $\overline{\sigma_i}\cap\overline{\sigma_j}=\varnothing$ для
всех $i\ne j$, $i,j= 1,2,\ldots$;

\medskip

(ii) $\sigma_{m-1}$ и $\sigma_{m+1}$ содержатся в различных
компонентах $D\setminus \sigma_m$ для всех $m>1$;

\medskip

(iii) $\cap\,d_m=\varnothing$, где $d_m$ -- компонента $D\setminus
\sigma_m$, содержащая $\sigma_{m+1}$.

\medskip
Согласно определению, цепь разрезов $\{\sigma_m\}$ определяет цепь
областей $d_m\subset D$, таких, что $\partial\,d_m\cap
D\subset\sigma_m$ и $d_1\supset d_2\supset\ldots\supset
d_m\supset\ldots$. Две цепи разрезов $\{\sigma_m\}$ и
$\{\sigma_k^{\,\prime}\}$ называются {\it эквивалентными}, если для
каждого $m=1,2,\ldots$ область $d_m$ содержит все области
$d_k^{\,\prime}$ за исключением конечного числа и для каждого
$k=1,2,\ldots$ область $d_k^{\,\prime}$ также содержит все области
$d_m$ за исключением конечного числа. {\it Конец} области $D$ -- это
класс эквивалентных цепей разрезов $D$.

Пусть $K$ -- конец области $D$ в ${\Bbb R}^n$, $\{\sigma_m\}$ и
$\{\sigma_m^{\,\prime}\}$ -- две цепи в $K$, $d_m$ и
$d_m^{\,\prime}$ -- области, соответствующие $\sigma_m$ и
$\sigma_m^{\,\prime}$. Тогда
$$\bigcap\limits_{m=1}\limits^{\infty}\overline{d_m}\subset
\bigcap\limits_{m=1}\limits^{\infty}\overline{d_m^{\,\prime}}\subset
\bigcap\limits_{m=1}\limits^{\infty}\overline{d_m}\ ,$$ и, таким
образом,
$$\bigcap\limits_{m=1}\limits^{\infty}\overline{d_m}=
\bigcap\limits_{m=1}\limits^{\infty}\overline{d_m^{\,\prime}}\ ,$$
т.\,е. множество
$$I(K)=\bigcap\limits_{m=1}\limits^{\infty}\overline{d_m}$$ зависит
только от $K$ и не зависит от выбора цепи разрезов $\{\sigma_m\}$.
Множество $I(K)$ называется {\it телом конца} $K$.

\medskip
Борелева функция $\rho:{\Bbb R}^n\,\rightarrow [0,\infty]$
называется {\it допустимой} для семейства $\Gamma$ кривых $\gamma$ в
${\Bbb R}^n,$ если $\int\limits_{\gamma}\rho (x)|dx|\ge 1$ для всех
(локально спрямляемых) кривых $ \gamma \in \Gamma$ (т.е.,
произвольная кривая $\gamma$ семейства $\Gamma$ имеет длину, не
меньшую $1$ в метрике $\rho$). В этом случае мы пишем: $\rho \in
{\rm adm} \,\Gamma.$ Для фиксированного $p\geqslant 1$ {\it
$p$-модулем} семейства кривых $\Gamma $ называется величина
$$M_p(\Gamma)=\inf_{\rho \in \,{\rm adm}\,\Gamma}
\int\limits_D \rho ^p (x) dm(x)\,.$$
При этом, если ${\rm adm}\Gamma=\varnothing,$ то полагаем:
$M_p(\Gamma)=\infty.$ Положим, кроме того, $M(\Gamma);=M_n(\Gamma).$

Далее, как обычно, для множеств $A$, $B$ и $C$ в ${\Bbb R}^n$,
$\Gamma(A,B,C)$ обозначает семейство всех кривых, соединяющих $A$ и
$B$ в $C$.

Следуя \cite{Na}, будем говорить, что конец $K$ является {\it
простым концом}, если $K$ содержит цепь разрезов $\{\sigma_m\}$,
такую, что
\begin{equation}\label{eqSIMPLE}\lim\limits_{m\rightarrow\infty}
M(\Gamma(C, \sigma_m, D))=0
\end{equation}
для некоторого континуума $C$ в $D$, где $M$ -- модуль семейства
$\Gamma(C, \sigma_m, D).$

\medskip
Будем говорить, что граница области $D$ в ${\Bbb R}^n$ является {\it
локально квазиконформной}, если каждая точка $x_0\in\partial D$
имеет окрестность $U$, которая может быть отображена квазиконформным
отображением $\varphi$ на единичный шар ${\Bbb B}^n\subset{\Bbb
R}^n$ так, что $\varphi(\partial D\cap U)$ является пересечением
${\Bbb B}^n$ с координатной гиперплоскостью. Говорим, что
ограниченная область $D$ в ${\Bbb R}^n$ {\it регулярна}, если $D$
может быть квазиконформно отображена на область с локально
квазиконформной границей.

\medskip
\begin{remark}\label{rem2}
Как следует из теоремы 4.1 в \cite{Na}, при квазиконформных
отображениях $g$ области $D_0$ с локально квазиконформной границей
на область $D$ в ${\Bbb R}^n$, $n\geqslant2$, существует
естественное взаимно однозначное соответствие между точками
$\partial D_0$ и простыми концами области $D$ и, кроме того,
предельные множества $C(g,b)$, $b\in\partial D_0$, совпадают с телом
$I(P)$ соответствующих простых концов $P$ в $D$.

Если $\overline{D}_P$ является пополнением регулярной области $D$ ее
простыми концами и $g_0$ является квазиконформным отображением
области $D_0$ с локально квазиконформной границей на $D$, то оно
естественным образом определяет в $\overline{D}_p$ метрику
$\rho_0(p_1,p_2)=\left|{\widetilde
{g_0}}^{-1}(p_1)-{\widetilde{g_0}}^{-1}(p_2)\right|$, где
${\widetilde {g_0}}$ продолжение $g_0$ в $\overline {D_0}$,
упомянутое выше.

Если $g_*$ является другим квазиконформным отображением некоторой
области $D_*$ с локально квазиконформной границей на область $D$, то
соответствующая метрика
$\rho_*(p_1,p_2)=\left|{\widetilde{g_*}}^{-1}(p_1)-{\widetilde{g_*}}^{-1}(p_2)\right|$
порождает ту же самую сходимость и, следовательно, ту же самую
топологию в $\overline{D}_P$ как и метрика $\rho_0$, поскольку
$g_0\circ g_*^{-1}$ является квазиконформным отображением между
областями $D_*$ и $D_0$, которое по теореме 4.1 из \cite{Na}
продолжается до гомеоморфизма между $\overline{D_*}$ и
$\overline{D_0}$.
\end{remark}

В дальнейшем, будем называть данную топологию в пространстве
$\overline{D}_P$ {\it топологией простых концов} и понимать
непрерывность отображений
$F:\overline{D}_P\rightarrow\overline{D^{\,\prime}}_P$ как раз
относительно этой топологии.

\medskip
Отображение $f:D\rightarrow \overline{{\Bbb R}^n}$ условимся
называть {\it $Q$-отоб\-ра\-же\-ни\-ем}, если $f$ удовлетворяет
соотношению
\begin{equation}\label{eq2*!}
M(f(\Gamma))\le \int\limits_D Q(x)\cdot \rho^n (x) dm(x)
\end{equation}
для произвольного семейства кривых $\Gamma$ в области $D$ и каждой
допустимой функции $\rho\in {\rm adm\,}\Gamma.$

Пусть $(X,d)$ и $\left(X^{\,{\prime}},{d}^{\,{\prime}}\right)$ ~---
метрические пространства с расстояниями $d$  и ${d}^{\,{\prime}},$
соответственно. Семейство $\frak{F}$ отображений $f:X\rightarrow
{X}^{\,\prime}$ называется {\it равностепенно непрерывным в точке}
$x_0 \in X,$ если для любого $\varepsilon > 0$ найдётся $\delta
> 0,$ такое, что ${d}^{\,\prime} \left(f(x),f(x_0)\right)<\varepsilon$ для
всех $f \in \frak{F}$ и  для всех $x\in X$ таких, что
$d(x,x_0)<\delta.$ Говорят, что $\frak{F}$ {\it равностепенно
непрерывно}, если $\frak{F}$ равностепенно непрерывно в каждой точке
из $x_0\in X.$ Всюду далее, если не оговорено противное, $d$ -- одна
из метрик в пространстве простых концов относительно области $D,$
упомянутых выше, а $d^{\,\prime}$ -- евклидова метрика.

\medskip
Для областей $D,$ $D^{\,\prime}\subset {\Bbb R}^n,$ $z_0, x_0\in D,$
$z_0\ne x_0,$ $z_0^{\,\prime}, x_0^{\,\prime}\in D^{\,\prime},$
$z_0^{\prime}\ne x_0^{\prime},$ и произвольной измеримой по Лебегу
функции $Q(x): {\Bbb R}^n\rightarrow [0, \infty],$ $Q(x)\equiv 0$
при $x\not\in D,$ обозначим через ${\frak H}_{z_0, x_0,
z_0^{\,\prime}, x_0^{\,\prime}, Q }(D, D^{\,\prime})$ семейство всех
гомеоморфизмов $f:D\rightarrow D^{\,\prime},$ $f(D)=D^{\,\prime},$
удовлетворяющих соотношению (\ref{eq2*!}) для произвольной $\rho\in
{\rm adm\,}\Gamma,$ таких, что $f(z_0)=z_0^{\,\prime},$
$f(x_0)=x_0^{\,\prime}.$ Имеет место следующее утверждение.

\medskip
\begin{theorem}\label{th4} {\sl Предположим, что
область $D$ в ${\Bbb R}^n,$ $n\geqslant 3,$ регулярна, а область
$D^{\,\prime}$ ограничена, имеет локально квазиконформную границу и,
одновременно, является $QED$-об\-лас\-тью, $Q\in L^1(D)$ и, кроме
того, $I(P_1)\cap I(P_2)=\varnothing$ для любых различных простых
концов $P_1, P_2\subset E_D,$ где, как обычно, $I(P)$ обозначает
тело простого конца $P\subset E_D.$
Тогда каждый элемент $g$ семейства ${\frak H}^{\,-1}_{z_0, x_0,
z_0^{\,\prime}, x_0^{\,\prime}, Q }(D, D^{\,\prime}),$ состоящего из
всех обратных гомеоморфизмов
$\left\{g=f^{-1}:D^{\,\prime}\rightarrow D\right\},$ где отображение
$f\in {\frak H}_{z_0, x_0, z_0^{\,\prime}, x_0^{\,\prime}, Q }(D,
D^{\,\prime}),$ может быть продолжен по непрерывности до отображения
$\overline{g}=\overline{f^{\,-1}}:
\overline{D_P^{\,\prime}}\rightarrow\overline{D_P},$ причём
семейство $\overline{{\frak H}^{\,-1}_{z_0, x_0, z_0^{\,\prime},
x_0^{\,\prime}, Q }(D, D^{\,\prime})},$ состоящее из всех
продолженных, таким образом, отображений
$\overline{g}:\overline{D_P^{\,\prime}}\rightarrow \overline{D_P},$
$g\in {\frak H}^{\,-1}_{z_0, x_0, z_0^{\,\prime}, x_0^{\,\prime}, Q
}(D, D^{\,\prime}),$ является равностепенно непрерывным в
$\overline{D^{\,\prime}_P}.$}
\end{theorem}

\medskip
{\bf 2. Леммы о непрерывном продолжении обратных отображений.} Как
обычно, {\it жордановой дугой} будем называть гомеоморфизм
$\gamma:I\rightarrow {\Bbb R}^n,$ где $I$ -- отрезок в ${\Bbb R}.$
Следуя при $n=2$ доказательству леммы 3, $\S\,28$ в \cite{Am}, а
также следствие 1.5.IV в \cite{HW} при $n\geqslant 3,$ имеем
следующее утверждение.

\medskip
\begin{lemma}\label{lem1}
{\sl Жорданова дуга не разбивает область $D$ в ${\Bbb R}^n,$
$n\geqslant 2.$ Если $D_1$ -- подобласть области $D\subset {\Bbb
R}^n,$ $n\geqslant 3,$ и $\gamma:I\rightarrow D$ -- жорданова дуга,
то также $\gamma$ не разбивает и $D_1.$}
\end{lemma}

\medskip
Пусть $E,$ $F\subset \overline{{\Bbb R}^n}$ -- произвольные
множества. Обозначим через $\Gamma(E,F,D)$ семейство всех кривых
$\gamma:[a,b]\rightarrow\overline{{\Bbb R}^n},$ которые соединяют
$E$ и $F$ в $D,$ т.е. $\gamma(a)\in E,$ $\gamma(b)\in F$ и
$\gamma(t)\in D$ при $t\in (a, b).$ Дальнейшее изложение существенно
опираются на аппарат нак называемых кольцевых $Q$-гомеоморфизмов
(см. \cite[глава~7]{MRSY}). Дадим определение этого класса
отображений. Пусть $x_0\in \overline{D},$ тогда отображение
$f:D\rightarrow \overline{{\Bbb R}^n}$ будем называть {\it кольцевым
$Q$-отображением относительно в точке $x_0$,}  для некоторого
$r_0=r(x_0)$ и произвольных сферического кольца $A=A(r_1,r_2,x_0),$
центрированного в точке $x_0,$ радиусов: $r_1,$ $r_2,$ $0<r_1<r_2<
r_0=r(x_0)$ и любых континуумов $E_1\subset \overline{B(x_0,
r_1)}\cap D,$ $E_2\subset \left(\overline{{\Bbb R}^n}\setminus
B(x_0, r_2)\right)\cap D$ отображение $f$ удовлетворяет соотношению
\begin{equation}\label{eq3*!!}
 M\left(f\left(\Gamma\left(E_1,\,E_2,\,D\right)\right)\right)\ \le
\int\limits_{A} Q(x)\cdot \eta^n(|x-x_0|)\ dm(x) \end{equation}
для каждой измеримой функции $\eta :(r_1,r_2)\rightarrow [0,\infty
],$ такой что
\begin{equation}\label{eq28*}
\int\limits_{r_1}^{r_2}\eta(r)dr \ge\ 1\,.
\end{equation}
Следующее утверждение для внутренних точек области $x_0$ области $D$
доказано в \cite{MRSY}, см. лемму 7.3, а для точек $x_0\in \partial
D$ может быть установлено по аналогии.

\medskip
\begin{proposition}\label{pr1} {\sl Пусть $D$ -- область в ${\Bbb R}^n,$ $Q:D\rightarrow
[0\,,\infty]$ -- измеримая функция, $q_{x_0}(r)$
--- среднее значение $Q(x)$ над сферой $|x-x_0|=r.$ Положим
$$
I=I(x_0,r_1,r_2)=\int\limits_{r_1}^{r_2}
\frac{dr}{rq_{x_0}^{\frac{1}{n-1}}(r)}
$$
и $S_j=\{x\in{\Bbb R}^n: |x-x_0|=r_j\},$ $j=1,2\,,$ где $x_0\in
\overline{D}$. Тогда для любого кольцевого
$Q$-го\-ме\-о\-мор\-физ\-ма $f:D\rightarrow\overline{{\Bbb R}^n}$ в
точке $x_0$
$$
M\left(\Gamma\left(f(S_1),f(S_2), f(D)\right)\right)\leqslant
\frac{\omega_{n-1}}{I^{n-1}}\,,$$
где $\omega_{n-1}$ --- площадь единичной сферы в ${\Bbb R}^n.$}
\end{proposition}

\medskip
Следующее утверждение сформулировано и доказано в \cite{KR}, см.
лемму 2.

\medskip
\begin{lemma}\label{thabc3} {\sl Каждый простой конец $P$ регулярной
области $D$ в ${\Bbb R}^n$, $n\geqslant2$, содержит цепь разрезов
$\sigma_m$, лежащую на сферах $S_m$ с центром в точке
$x_0\in\partial D$ и с евклидовыми радиусами $r_m\rightarrow 0$ при
$m\rightarrow\infty$.}\end{lemma}

\medskip
Следуя \cite{KR}, говорим, что граница $\partial D$ является {\it
слабо плоской в точке} $x_0\in\partial D$, если для каждой
окрестности $U$ точки $x_0$ и для любого числа $P>0$ найдется
окрестность $V\subset U$ точки $x_0$, такая, что
\begin{equation}\label{eq1.5KR} M(\Gamma(E, F, D))\geqslant P\end{equation}
для любых континуумов $E$ и $F$ в $D$, пересекающих $\partial U$ и
$\partial V$. Говорим также, что граница $\partial D$ является {\it
слабо плоской}, если она слабо плоская в каждой точке из $\partial
D$.

\medskip
\begin{remark}\label{rem1}
Отметим, что произвольные $QED$-об\-лас\-ти в ${\Bbb R}^n$ имеют
слабо плоские границы (см., напр., \cite[замечание~13.10 и
разд.~13.10]{MRSY}).
\end{remark}

\medskip
Следующие два утверждения доказаны в работе \cite{KR} в пространстве
${\Bbb R}^n$ для несколько иных классов отображений (см. лемму 4 и
теорему 1), а также на плоскости в случае тех же классов (см. лемму
6.1 и теорему 6.1 в \cite{GRY}). Доказательство этих утверждений не
содержит существенных отличий упомянутых случаев, однако, ради
полноты изложения мы приведём их полностью.

\medskip
\begin{lemma}\label{l:9.1} {\sl Пусть $D$ и $D^{\,\prime}$ -- регулярные области в
${\Bbb R}^n$, $n\geqslant 2$, $P_1$ и $P_2$ -- разные простые концы
области $D$, и пусть $\sigma_m$, $m=1,2,\ldots$, -- цепь разрезов
простого конца $P_1$ из леммы \ref{thabc3}, лежащая на сферах
$S(z_1,r_m)$, $z_1\in I(P_1)$, с ассоциированными областями $D_m$.
Предположим, что функция $Q$ интегрируема по сферам
\begin{equation}\label{INTERSECTION}
D(r)=\left\{x\in D:|x-z_1|=r\right\}=D\cap S(z_1,r)\end{equation}
для некоторого множества $E$ чисел $r\in(0,d)$ положительной
линейной меры, где $d=r_{m_0}$ и $m_0$ минимальное из чисел, таких,
что область $D_{m_0}$ не содержит последовательностей точек,
сходящихся к $P_2$. Если $f$ -- кольцевой $Q$-гомеоморфизм области
$D$ на область $D^{\,\prime}$ в точке $z_1$ и $\partial
D^{\,\prime}$ является слабо плоской, то
\begin{equation}\label{e:9.2}C(P_1,f)\cap C(P_2,f)=\varnothing.\end{equation}}\end{lemma}

Заметим, что в силу метризуемости пополнения $\overline{D}_P$
области $D$ простыми концами, см. сделанное выше замечание, число
$m_0$ в лемме \ref{l:9.1} всегда существует.

\medskip
\begin{proof} Выберем $\varepsilon\in(0,d)$ так, что $E_0:=\{r\in
E:r\in(\varepsilon,d)\}$ имеет положительную линейную меру. Такой
выбор возможен в силу полуаддитивности линейной меры и исчерпания
$E=\cup E_m$, где $E_m=\{r\in E:r\in(1/m,d)\}$, $m=1,2,\ldots$.
Полагая $S_1=S(z_1, \varepsilon),$ $S_2=S(z_1, d),$ ввиду
предложения \ref{pr1} будем иметь:
\begin{equation}\label{e:9.3}
M(\Gamma(f(S_1), f(S_2), f(D)))<\infty\,.
\end{equation}

Предположим, что $C_1\cap C_2\ne\varnothing$, где $C_i=C(P_i,f)$,
$i=1,2$. По построению существует $m_1>m_0$, такое, что
$\sigma_{m_1}$ лежит на сфере $S(z_1,r_{m_1})$ с
$r_{m_1}<\varepsilon$. Пусть $D_0=D_{m_1}$ и $D_*\subseteq
D\setminus D_{m_0}$ -- область, ассоциированная с цепью разрезов
простого конца $P_2$. Пусть $y_0\in C_1\cap C_2$. Выберем $r_0>0$
так, что $S(y_0,r_0)\cap f(D_0)\ne\varnothing$ и $S(y_0,r_0)\cap
f(D_*)\ne\varnothing$.

Обозначим $\Gamma=\Gamma(\overline{D_0},\overline{D_*};D)$. Тогда по
принципу минорирования и из (\ref{e:9.3}) следует, что
\begin{equation}\label{e:9.4}
M(f(\Gamma))\leqslant M(\Gamma(f(S_1), f(S_2),
f(D)))<\infty\,.\end{equation}
Пусть $M_0>M(f(\Gamma))$ -- конечное число. Из условия слабой
плоскости $\partial D^{\,\prime}$ следует, что найдется
$r_*\in(0,r_0)$, такое, что
$$M(\Gamma(E, F, D^{\,\prime}))\geqslant M_0$$
для всех континуумов $E$ и $F$ в $D^{\,\prime}$, пересекающих сферы
$S(y_0,r_0)$ и $S(y_0,r_*)$. Однако, эти сферы могут быть соединены
непрерывными кривыми $c_1$ и $c_2$ в областях $f(D_0)$ и $f(D_*)$ и,
в частности, для этих кривых
\begin{equation}\label{e:9.4a}
M_0\leqslant M(\Gamma(c_1,c_2;D^{\,\prime}))\leqslant
M(f(\Gamma)).\end{equation} Полученное противоречие опровергает
предположение, что $C_1\cap C_2\ne\varnothing$. \end{proof}

\medskip
\begin{theorem}\label{t:9.5} {\it Пусть $D$ и $D^{\,\prime}$ -- регулярные области в ${\Bbb R}^n$,
$n\geqslant 2$. Если $f$ -- кольцевой $Q$-гомеоморфизм $D$ на
$D^{\,\prime}$ в каждой точке $x_0\in \partial D$ и $Q\in L(D)$, то
$f^{\,-1}$ продолжается до непрерывного отображения
$\overline{D^{\,\prime}}_P$ на $\overline{D}_P$.}
\end{theorem}

\begin{proof} По теореме Фубини, см., напр., \cite{Sa}, множество
$$E(x_0)=\left\{r\in(0,d(x_0)):Q|_{D(x_0,r)}\in
L(D(x_0,r))\right\} \quad\forall\ x_0\in\partial D,$$ где
$d(x_0)=\sup_{x\in D}|x-x_0|$ и $D(x_0,r)=D\cap S(x_0,r)$, имеет
положительную линейную меру, поскольку $Q\in L(D)$. Согласно
сделанном во введении замечаниям,  без ограничения общности можем
считать, что область $D^{\,\prime}$ имеет слабо плоскую границу.
Таким образом, рассуждая от противного и с учетом метризуемости
пространств $\overline{D^{\,\prime}}_P$ и $\overline{D}_P$ ,
получаем заключение о продолжении $f^{\,-1}$ на
$\overline{D^{\,\prime}}.$

Таким образом, мы имеем продолжение $f^{\,-1}$ на
$\overline{D^{\,\prime}}$ такое, что $C(\partial D^{\,\prime},
f^{\,-1})\subset \overline{D}_P\setminus D.$ Покажем, что
$C(\partial D^{\,\prime}, f^{\,-1})=\overline{D}_P\setminus D.$
Действительно, если $P_0$ -- простой конец в $D,$ то  найдётся
последовательность $x_n\rightarrow P_0$ при $n\rightarrow \infty.$
Ввиду компактности $\overline{D}$ и $\overline{D^{\,\prime}}$ мы
можем считать, что $x_n\rightarrow x_0\in \partial D$ и
$f(x_n)\rightarrow \zeta_0\in \partial D^{\,\prime}$ при
$n\rightarrow\infty.$ Последнее означает, что $P_0\in C(\zeta_0,
f^{\,-1}).$

Покажем, наконец, что продолженное отображение
$g:\overline{D^{\,\prime}}\rightarrow \overline{D}_P$ непрерывно.
Действительно, пусть $\zeta_n\rightarrow \zeta_0$ при
$n\rightarrow\infty,$ $\zeta_n, \zeta_0\in \overline{D^{\,\prime}}.$
Если $\zeta_0$ -- внутренняя точка области $D^{\,\prime},$ желанное
утверждение очевидно. Пусть $\zeta_0\in \partial D^{\,\prime},$
тогда выберем $\zeta_n^*\in D^{\,\prime}$ таким, что
$|\zeta_n-\zeta_n^*|<1/n$ и $\rho(g(\zeta_n), g(\zeta_n^*))<1/n,$
где $\rho$ -- одна из метрик, указанных в замечании \ref{rem2}. По
построению $g(\zeta_n^*)\rightarrow g(\zeta_0),$ поскольку
$\zeta_n^*\rightarrow \zeta_0$ и, значит, также и
$g(\zeta_n)\rightarrow g(\zeta_0)$ при $n\rightarrow\infty.$~$\Box$
\end{proof}

\medskip
{\it Доказательство теоремы \ref{th4}.} {\bf I. } Не ограничивая
общности рассуждений, ввиду замечания \ref{rem2} можно считать, что
$\overline{D^{\,\prime}_P}=\overline{D^{\,\prime}}.$

Т.к. область $D^{\,\prime}$ является $QED$-об\-лас\-тью, каждый
обратный гомеоморфизм $f^{\,-1}$ имеет непрерывное продолжение на
границу $D^{\,\prime}$ ввиду теоремы \ref{t:9.5}. Осталось показать
равностепенную непрерывность семейства отображений ${\frak
H}^{\,-1}_{z_0, x_0, z_0^{\,\prime}, x_0^{\,\prime}, Q }(D,
D^{\,\prime})$ (обозначения не меняем) в области
$\overline{D^{\,\prime}}.$

{\bf II. } Покажем сначала, что  ${\frak H}^{\,-1}_{z_0, x_0,
z_0^{\,\prime}, x_0^{\,\prime}, Q }(D, D^{\,\prime})$ равностепенно
непрерывно в $\overline{D^{\,\prime}}\setminus \{z_0^{\,\prime},
x_0^{\,\prime}\}.$ Предположим противное, т.е. найдётся $y_0\in
\overline{D^{\,\prime}},$ $x_0^{\,\prime}\ne y_0\ne z_0^{\,\prime},$
и $\varepsilon_0>0,$ такие что для любого $m\in {\Bbb N}$ существует
$y_m\in \overline{D^{\,\prime}}$ с $|y_m-y_0|<1/m$ и элемент
$f_m^{\,-1}\in {\frak H}^{\,-1}_{z_0, x_0, z_0^{\,\prime},
x_0^{\,\prime}, Q }(D, D^{\,\prime}),$  такие что
\begin{equation}\label{eq13***}
\rho(f_m^{\,-1}(y_m), f_m^{\,-1}(y_0))\ge \varepsilon_0\,,
\end{equation}
где $\rho$ -- одна из метрик из замечания \ref{rem2}. Так как
$f_m^{\,-1}$ непрерывным образом продолжаются на
$\overline{D^{\,\prime}},$ найдутся последовательности $z_m, x_m\in
D^{\,\prime}$ такие, что
$$|z_m-y_m|<1/m\,,\quad |x_m-y_0|<1/m\,, \quad m\rightarrow\infty\,,$$
и, при этом,
\begin{equation}\label{eq4}
\rho(f_m^{\,-1}(y_m), f_m^{\,-1}(z_m))<1/m\,,\quad
\rho(f_m^{\,-1}(x_m), f_m^{\,-1}(y_0))<1/m\,.
\end{equation}
Тогда из (\ref{eq13***}) вытекает, что
\begin{equation}\label{eq13}
\rho(f_m^{\,-1}(z_m), f_m^{\,-1}(x_m))\ge \varepsilon_1\,,
\end{equation}
где $\varepsilon_1>0$ -- некоторое фиксированное число.
Т.к. $\overline{D}_P$ является компактом, то можно считать, что для
некоторых $P_0^1, P_0^2\in \overline{D}_P$ выполнены условия
\begin{equation}\label{eq5}
f_m^{\,-1}(z_m)\rightarrow P_0^1\,,\quad f_m^{\,-1}(x_m)\rightarrow
P_0^2\,,\quad m\rightarrow\infty\,.
\end{equation}
В частности, из (\ref{eq4}) и (\ref{eq5}) в силу неравенства
треугольника следует, что
\begin{equation}\label{eq6}
f_m^{-1}(y_0)\rightarrow P_0^2\in \overline{D}_P\,, \quad
m\rightarrow\infty\,.
\end{equation}
В силу неравенства (\ref{eq13})
%
$$\rho(P_0^1, P_0^2)\ge \varepsilon_1/2\,.$$
%
Не ограничивая общности можно считать, что при всех $m\in {\Bbb N}$
выполнены включения $f_m^{\,-1}(z_m)\in D_m$ и $f_m^{\,-1}(x_m)\in
D_m^{\,\prime},$ где $D_m$ и $D_m^{\,\prime}$ -- последовательности
областей, соответствующие $P_0^1, P_0^2\in \overline{D}_P$ (в
случае, если $P_0^1$ либо $P_0^2$ -- внутренние точки области $D,$ в
качестве областей $D_m$ либо $D_m^{\,\prime}$ берём
последовательности открытых шаров, стягивающихся к $P_0^1$ либо
$P_0^2,$ соответственно). Можно считать, что $x_0\ne P_0^1$ и
$z_0\ne P_0^2,$ и что $D_k\cap D_l^{\,\prime}=\varnothing$ при всех
$k, l\in {\Bbb N}.$ Более того, поскольку по условию теоремы тела
простых концов $P_0^1$ и $P_0^2$ не пересекаются, то можно считать,
что
\begin{equation}\label{eq19}
\overline{D_i}\cap\overline{D_i^{\,\prime}}=\varnothing\quad\forall\,\,i\in
{\Bbb N}\,.
\end{equation}
Заметим, что (\ref{eq19}) верно также и в случае, когда в качестве
одного из $P_0^i$ берётся внутренняя точка области $D.$

\medskip
{\bf III. } Построим последовательность кривых $\alpha_m$ следующим
образом. Кривую $\alpha_1$ определим как произвольную дугу,
соединяющую точки $x_0$ и $f_1^{\,-1}(z_1)\in D_1$ в области
$D\setminus \overline{D_1^{\,\prime}}.$ Далее, соединим точки
$f_1^{\,-1}(z_1)$ и $f_2^{\,-1}(z_2)$ внутри области $D_1$ некоторой
дугой $\gamma_1,$ и определим $\alpha_2$ как объединение $\alpha_1$
и $\gamma_1.$ И так далее. На некотором $m$-м шаге построим кривую
$\alpha_m,$ которая будет определяться как кривая $\alpha_{m-1},$
объединённая с дугой $\gamma_{m-1},$ где $\gamma_{m-1}$ --
произвольная дуга, соединяющая точки $f_{m-1}^{\,-1}(z_{m-1})$ и
$f_m^{\,-1}(z_m)$ в области $D_{m-1}.$ И так далее.

\medskip
{\bf IV.} Построим также последовательность кривых $\beta_m$
следующим образом. Кривую $\beta_1$ определим как произвольную дугу,
соединяющую точки $z_0$ и $f_1^{\,-1}(x_1)\in D^{\prime}_1$ в
области $D\setminus \overline{D_1}$  так, чтобы
$|\alpha_1|\cap|\beta_1|=\varnothing,$ что возможно ввиду леммы
\ref{lem1}. Далее, соединим точки $f_1^{\,-1}(x_1)$ и
$f_2^{\,-1}(x_2)$ внутри области $D^{\,\prime}_1$ некоторой дугой
$\delta_1,$ и определим $\beta_2$ как объединение  $\beta_1$ и
$\delta_1.$ И так далее. На некотором $m$-м шаге построим кривую
$\beta_m,$ которая будет определяться как кривая $\beta_{m-1},$
объединённая с дугой $\delta_{m-1},$ где $\delta_{m-1}$ --
произвольная дуга, соединяющая точки $f_{m-1}^{\,-1}(x_{m-1})$ и
$f_m^{\,-1}(x_m)$ в области $D^{\,\prime}_{m-1}.$ И так далее.

\medskip
{\bf V.} Заметим, что найдётся постоянная $C>0$ такая, что
\medskip
\begin{equation}\label{eq1}
{\rm dist}(|\alpha_m|, |\beta_m|)\geqslant C\qquad \forall\quad m\in
{\Bbb N}\,,
\end{equation}
где $$|\gamma|:=\{x\in D:\exists\,t: \gamma(t)=x\}\,.$$

Действительно, по построению $|\alpha_m|\subset |\alpha_1|\cup
\overline{D_1},$ а $|\beta_m|\subset |\beta_1|\cup
\overline{D^{\,\prime}_1},$ причём ввиду соотношения (\ref{eq19})
множества $C_1:=|\alpha_1|\cup \overline{D_1}$ и $C_2:=|\beta_1|\cup
\overline{D^{\,\prime}_1}$ являются непересекающимися компактными
подмножествами в ${\Bbb R}^n$ а, значит, отстоят друг от друга на
расстояние не меньшее некоторого $C>0.$ Тем более, $|\alpha_m|$ и
$|\beta_m|$ отстоят друг от друга не меньше, чем на $C,$ что и
доказывает соотношение (\ref{eq1}).

\medskip
{\bf VI.}  Пусть $\Gamma_m$ ~--- семейство кривых, соединяющих
множества $\alpha_m$ и $\beta_m$ в $D,$ тогда функция
$$\rho(x)= \left\{
\begin{array}{rr}
\frac{1}{C}, & x\in D\\
0,  &  x\notin  D
\end{array}
\right. $$
является допустимой для семейства $\Gamma_m,$ и кроме того, т.к.
$f_m$ являются $Q$-го\-ме\-о\-мор\-физ\-ма\-ми в $D,$
\begin{equation}\label{eq14***}
M(f_m(\Gamma_m))\le \frac{1}{C^n}\int\limits_{D}
Q(x)\,dm(x):=c(C)<\infty\,,
\end{equation}
т.к. $Q\in L^1(D).$

\medskip
{\bf VII.} Будем иметь: $z_m, x^{\,\prime}_0\in f_m(\alpha_m)$ и,
поскольку по построению, $z_m\rightarrow y_0\ne x_0^{\,\prime},$ то
найдётся постоянная $\delta_1>0$ такая, что ${\rm
diam\,}(f_m(\alpha_m))\geqslant \delta_1>0$ при всех $m\in {\Bbb
N}.$ Аналогично, $x_m, z^{\,\prime}_0\in f_m(\beta_m)$ и, поскольку
по построению, $x_m\rightarrow y_0\ne z_0^{\,\prime},$ то найдётся
постоянная $\delta_2>0$ такая, что ${\rm
diam\,}(f_m(\beta_m))\geqslant \delta_2>0$ при всех $m\in {\Bbb N}.$
Кроме того, заметим, что ${\rm dist\,}(f_m(\alpha_m),
f_m(\beta_m))\rightarrow 0$ при $m\rightarrow\infty.$ Согласно
свойству сближающихся континуумов (см. теорему 2.3 и замечание 2.8 в
\cite{MRSY}) будем иметь, что
\begin{equation}\label{eq2}
M(\Gamma(f_m(\alpha_m),f_m(\beta_m), {\Bbb R}^n))\rightarrow\infty,
\quad m\rightarrow\infty\,.
\end{equation}
Поскольку $D^{\,\prime}$ является $QED$-областью, то из (\ref{eq2})
вытекает, что $M(f_m(\Gamma_m))\rightarrow\infty$ при
$m\rightarrow\infty,$ что противоречит соотношению (\ref{eq14***}).
Полученное противоречие говорит о том, что семейство отображений
${\frak H}^{\,-1}_{z_0, x_0, z_0^{\,\prime}, x_0^{\,\prime}, Q }(D,
D^{\,\prime})$ равностепенно непрерывно в точке $y_0.$ $\Box$

\medskip
{\bf VIII.} Для завершения доказательства нам осталось показать, что
семейство отображений ${\frak H}^{\,-1}_{z_0, x_0, z_0^{\,\prime},
x_0^{\,\prime}, Q }(D, D^{\,\prime})$ также является равностепенно
непрерывным в точках $x_0^{\,\prime}$ и $z^{\,\prime}_0.$ Рассмотрим
для определённости случай точки $x^{\,\prime}_0$ (случай точки
$z^{\,\prime}_0$ рассматривается аналогично). Предположим противное,
т.е. найдётся $\varepsilon_0>0$ такое, что для любого $m\in {\Bbb
N}$ существует $x_m\in D^{\,\prime}$ с $|x_m-x^{\,\prime}_0|<1/m$ и
элемент $f_m^{\,-1}\in {\frak H}^{\,-1}_{z_0, x_0, z_0^{\,\prime},
x_0^{\,\prime}, Q }(D, D^{\,\prime}),$  такие что
\begin{equation}\label{eq7}
|f_m^{\,-1}(x_m)- x_0|\ge \varepsilon_0\,.
\end{equation}
Поскольку пространство $\overline{D}_P$ компактно, мы можем считать,
что для некоторого $P_0\in \overline{D}_P,$ $P_0\ne x_0$ выполнено
\begin{equation}\label{eq8*}
\rho(f_m^{\,-1}(x_m), P_0)\rightarrow 0\,, \quad
m\rightarrow\infty\,.
\end{equation}
Пусть $D_m,$ $m=1,2,\ldots,$ -- последовательность областей,
соответствующих простому концу $P_0$ (если это точка в области $D,$
то, как и прежде, в качестве последовательности $D_m$ берём
сферические окрестности, сжимающиеся в точку). Можно считать, что
$f_m^{\,-1}(x_m)\in D_m$ для всех $m\in {\Bbb N}.$

\medskip
Поскольку по условию область $D$ регулярна, то она квазиконформно
отображается на некоторую область $D_0$ с локально квазиконформной
границей, которая по теореме Лиувилля не может совпадать с ${\Bbb
R}^n,$ а также, ввиду локальной квазиконформности границы, с ${\Bbb
R}^n\setminus\{b\}$ для некоторой точки $b\in {\Bbb R}^n.$ Ввиду
наличия взаимно однозначного соответствия между $\partial D_0$ и
$E_D,$ найдутся не менее двух простых концов $P_1, P_2\subset E_D,$
$P_1\ne P_2.$ Выберем в качестве $P_1\subset E_D$ произвольный
простой конец, не совпадающий с $P_0.$ В качестве вспомогательной
последовательности $y_m$ рассмотрим произвольную последовательность,
сходящуюся к $P_1.$

\medskip
В силу компактности $\overline{D_P}$ можно считать, что
последовательность $y_m$ сходится к некоторой граничной точке
$\zeta_0\in \partial D.$ Поскольку $f_1$ -- гомеоморфизм, то $C(f_1,
\zeta_0)\subset \partial D^{\,\prime}.$ Ввиду компактности
$\overline{D^{\,\prime}}$ найдётся подпоследовательность номеров
$m^1_k,$ такая что $f_1(y_{m^1_k})$ сходится к некоторой граничной
точке $\xi_1\in \partial D^{\,\prime}.$ Тогда найдётся номер $k_1\in
{\Bbb N}$ такой, что ${\rm dist}\, (f_1(y_{m^1_{k_1}}),
\partial D^{\,\prime})<1.$ Полагаем $l_1:=m^1_{k_1}.$ Рассмотрим последовательность $y_m,$
$m> l_1.$ Поскольку $f_2$ -- гомеоморфизм, то $C(f_2,
\zeta_0)\subset
\partial D^{\,\prime}.$ Ввиду компактности $\overline{D^{\,\prime}}$
найдётся подпоследовательность номеров $m^2_k,$ такая что
$f_2(y_{m^2_k})$ сходится к некоторой граничной точке $\xi_2\in
\partial D^{\,\prime}.$ Тогда найдётся номер $k_2\in
{\Bbb N}$ такой, что ${\rm dist}\, (f_2(y_{m^2_{k_2}}),
\partial D^{\,\prime})<1/2.$ Полагаем $l_2:=m^2_{k_2}.$ И так далее.
В результате бесконечного процесса получаем последовательность
$y_{l_k}$ такую, что ${\rm dist}\, (f_k(y_{l_k}),
\partial D^{\,\prime})<1/k.$ Полагаем $z_k:=y_{l_k}.$ Тогда
\begin{equation}\label{eq9}
{\rm dist}\, (f_k(z_k),
\partial D^{\,\prime})<1/k\,,k=1,2,\ldots,
\end{equation}
причём последовательность $z_k$ также сходится к простому концу
$P_1.$ Переходя, если нужно, к подпоследовательности, мы можем
считать, что $z_i\in D_i^{\,\prime},$ и что выполнено условие
(\ref{eq19}).

\medskip
Построим, как и в пунктах {\bf III} и {\bf IV}, последовательности
кривых $\alpha_m$ и $\beta_m$ следующим образом. Кривую $\alpha_1$
определим как произвольную дугу, соединяющую точки $x_0$ и $z_1$ в
$D\setminus \overline{D_1}.$ Далее соединим точки $z_1$ и $z_2$
внутри области $D^{\,\prime}_1$ некоторой дугой $\gamma_1,$ и
определим $\alpha_2$ как объединение $\alpha_1$ и $\gamma_1.$ И так
далее. На некотором $m$-м шаге построим кривую $\alpha_m,$ которая
будет определяться как кривая $\alpha_{m-1},$ объединённая с дугой
$\gamma_{m-1},$ где $\gamma_{m-1}$ -- произвольная дуга, соединяющая
точки $z_{m-1}$ и $z_m$ в области $D^{\,\prime}_{m-1}.$ И так далее.

\medskip
Построим также последовательность кривых $\beta_m$ следующим
образом. Кривую $\beta_1$ определим как произвольную дугу,
соединяющую точки $z_0$ и $f_1^{\,-1}(x_1)\in D^{\prime}_1$ в
области $D\setminus \overline{D^{\,\prime}_1}$  так, чтобы
$|\alpha_1|\cap|\beta_1|=\varnothing,$ что возможно ввиду леммы
\ref{lem1}. Далее, соединим точки $f_1^{\,-1}(x_1)$ и
$f_2^{\,-1}(x_2)$ внутри области $D_1$ некоторой дугой $\delta_1,$ и
определим $\beta_2$ как объединение  $\beta_1$ и $\delta_1.$ И так
далее. На некотором $m$-м шаге построим кривую $\beta_m,$ которая
будет определяться как кривая $\beta_{m-1},$ объединённая с дугой
$\delta_{m-1},$ где $\delta_{m-1}$ -- произвольная дуга, соединяющая
точки $f_{m-1}^{\,-1}(x_{m-1})$ и $f_m^{\,-1}(x_m)$ в области
$D_{m-1}.$ И так далее.

\medskip
Пусть $\Gamma_m$ ~--- семейство кривых, соединяющих множества
$\alpha_m$ и $\beta_m$ в $D.$ Аналогично тому, как это сделано в
пункте {\bf V}, доказывается справедливость соотношения вида
(\ref{eq1}), откуда вытекает оценка вида (\ref{eq14***}).

\medskip
Заметим, что в таком случае $f_k(z_k), x^{\,\prime}_0\in
f_k(\alpha_k)$ и, ввиду (\ref{eq9}) найдётся постоянная $\delta_1>0$
такая, что ${\rm diam\,}(f_k(\alpha_k))\geqslant \delta_1>0$ при
всех $k\geqslant k_0.$ Аналогично, $x_k, z^{\,\prime}_0\in
f_k(\beta_k)$ и, поскольку по построению, $x_k\rightarrow
x^{\,\prime}_0\ne z_0^{\,\prime},$ то найдётся постоянная
$\delta_2>0$ такая, что ${\rm diam\,}(f_k(\beta_k))\geqslant
\delta_2>0$ при всех $k\in {\Bbb N}.$ Кроме того, заметим, что ${\rm
dist\,}(f_k(\alpha_k), f_k(\beta_k))\rightarrow 0$ при
$k\rightarrow\infty.$ Согласно свойству сближающихся континуумов
(см. теорему 2.3 и замечание 2.8 в \cite{MRSY}) будем иметь, что
\begin{equation}\label{eq3}
M(\Gamma(f_k(\alpha_k),f_k(\beta_k), {\Bbb R}^n))\rightarrow\infty,
\quad k\rightarrow\infty\,.
\end{equation}
Поскольку $D^{\,\prime}$ является $QED$-областью, то из (\ref{eq3})
вытекает, что $M(f_k(\Gamma_k))\rightarrow\infty$ при
$k\rightarrow\infty,$ что противоречит соотношению (\ref{eq14***}).
Полученное противоречие говорит о том, что семейство отображений
${\frak H}^{\,-1}_{z_0, x_0, z_0^{\,\prime}, x_0^{\,\prime}, Q }(D,
D^{\,\prime})$ равностепенно непрерывно в точке $x_0.$ Теорема
доказана.~$\Box$

\medskip
КОНТАКТНАЯ ИНФОРМАЦИЯ

\medskip
\noindent{{\bf Руслан Радикович Салимов} \\
Институт математики НАН Украины \\
ул. Терещенковская, д. 3 \\
г. Киев-4, Украина, 01 601\\
тел. +38 095 630 85 92 (моб.), e-mail: ruslan623@yandex.ru}

\medskip
\noindent{{\bf Евгений Александрович Севостьянов} \\
Житомирский государственный университет им.\ И.~Франко\\
ул. Большая Бердичевская, 40 \\
г.~Житомир, Украина, 10 008 \\ тел. +38 066 959 50 34 (моб.),
e-mail: esevostyanov2009@mail.ru}

\end{document}